# A distribution for a pair of unit vectors generated by Brownian motion

SHOGO KATO

*Institute of Statistical Mathematics, 4-6-7 Minami-Azabu, Minato-ku, Tokyo 106-8569, Japan.*
*E-mail: skato@ism.ac.jp*

We propose a bivariate model for a pair of dependent unit vectors which is generated by Brownian motion. Both marginals have uniform distributions on the sphere, while the conditionals follow so-called "exit" distributions. Some properties of the proposed model, including parameter estimation and a pivotal statistic, are investigated. Further study is undertaken for the bivariate circular case by transforming variables and parameters into the form of complex numbers. Some desirable properties, such as a multiplicative property and infinite divisibility, hold for this submodel. Two estimators for the parameter of the submodel are studied and a simulation study is carried out to investigate the finite sample performance of the estimators. In an attempt to produce more flexible models, some methods to generalize the proposed model are discussed. One of the generalized models is applied to wind direction data. Finally, we show how it is possible to construct distributions in the plane and on the cylinder by applying bilinear fractional transformations to the proposed bivariate circular model.

## 1. Introduction

In a variety of scientific fields, observations are described as pairs of $d$-dimensional unit vectors. In meteorology, for example, wind directions at the weather station in Milwaukee at 6 a.m. and noon (Johnson and Wehrly (1977)) are data of this type with $d = 2$. Another example with $d = 3$ consists of directions of the magnetic field in a rock sample before and after some laboratory treatment (Stephens (1979)).

For the analysis of data of this type, various stochastic models have been proposed in the literature. Some distributions with certain marginals or conditionals are seen in Mardia (1975), Wehrly and Johnson (1980), Rivest (1988, 1997), Downs and Mardia (2002), SenGupta (2004) and Shieh and Johnson (2005). Among various works on distributions for bivariate angular data considerable attention has been paid to models with uniform marginals, which can be viewed as spherical equivalents of copulas. Johnson and Wehrly (1977) provided a bivariate circular distribution with uniform marginals and von Mises conditionals. Saw (1983) introduced a bivariate family with uniform marginals for pairs of dependent unit vectors which is an offset distribution of the multivariate normal distribution with some restrictions on parameters. Rivest (1984) discussed a certain class of distributions with so-called "$O(d)$-symmetric" densities, which has uniform marginals. Recent work by Alfonsi and Brigo (2005) proposed new families of copulas based on periodic functions.

The potential application of these special copulas is not restricted to the bivariate angular data whose marginals are uniformly distributed. Saw (1983) constructed a method which extends these copulas to distributions having any rotationally symmetric marginals. In the bivariate circular case, it is also possible to use a well-known technique in copula theory to generalize the model. These methods enable us to provide a bivariate model with prescribed marginals.

The main purpose of this paper is to introduce a new distribution with uniform marginals which is generated by $\mathbb{R}^d$-valued Brownian motion. To our knowledge, distributions for bivariate angular data have not previously been proposed based on Brownian motion. In this paper, a new approach is taken to provide a tractable model. This method enables us to define a class of bivariate distributions with uniform marginals and derive some desirable properties. We also discuss generalization of the proposed model.

Section 2 suggests a model for two dependent unit vectors and Section 3 investigates properties of the proposed model, including parameter estimation and a pivotal statistic. In Section 4, we focus on the bivariate circular case of the model and discuss its properties in detail. It is shown that some desirable properties, such as a multiplicative property and infinite divisibility, hold for this submodel. Some properties of two estimators for the parameter of the submodel are studied by means of a simulation study. In Section 5, generalizations of the proposed model are discussed and one of the extended models is applied to wind direction data. Finally, in Section 6, related models on $\mathbb{R}^2$ and on the cylinder are constructed by applying bilinear fractional transformations to the proposed model.

## 2. A model for a pair of unit vectors

### 2.1. Definition of the proposed model

Let $\{B_t; t \geq 0\}$ be $\mathbb{R}^d$-valued Brownian motion where $d \geq 2$. Starting at $B_0 = 0$, a Brownian particle will eventually hit a $d$-sphere with radius $\rho \, (\in (0, 1))$. Let $\tau_1$ be the minimum time at which the particle exits the sphere, that is, $\tau_1 = \inf\{t; \|B_t\| = \rho\}$, where $\|\cdot\|$ is the Euclidean norm. After leaving the sphere with radius $\rho$, the particle will hit a unit sphere first at time $\tau_2$, meaning $\tau_2 = \inf\{t; \|B_t\| = 1\}$. The proposed model is then defined by the joint distribution of a random vector

$$\left(Q\frac{B_{\tau_1}}{\|B_{\tau_1}\|}, B_{\tau_2}\right),$$

where $Q$ is a member of $O(d)$, the group of orthogonal transformations on $\mathbb{R}^d$. Note that the reason for multiplying $Q$ by $B_{\tau_1}/\|B_{\tau_1}\|$ is to make the model more flexible without losing its tractability.

## 2.2. Probability density function

For convenience, write $(U, V) = (QB_{\tau_1}/\|B_{\tau_1}\|, B_{\tau_2})$. It is clear that $(U, V)$ is a random vector where each variable takes values on the unit sphere. The joint distribution of $(U, V)$ has density

$$c(u, v) = \frac{1}{A_{d-1}^2} \frac{1 - \rho^2}{(1 - 2\rho u'Qv + \rho^2)^{d/2}}, \qquad u, v \in S^{d-1}, \tag{2.1}$$

where $\rho \in [0, 1)$, $Q \in O(d)$, $S^{d-1} = \{x \in \mathbb{R}^d; \|x\| = 1\}$, $x'$ is the transpose of $x$ and $A_{d-1}$ is surface area of $S^{d-1}$, that is, $A_{d-1} = 2\pi^{d/2}/\Gamma(d/2)$. The domain of $\rho$ is extended to include $\rho = 0$ so that the model includes the uniform distribution. We write $(U, V) \sim BS_d(\rho Q)$ if a random vector $(U, V)$ has density (2.1). For the derivation of the density (2.1), see Appendix A.

The parameter $\rho$ influences the dependence between $U$ and $V$. When $\rho = 0$, $U$ and $V$ are independent and distributed as the uniform distribution on the sphere, that is, $c(u, v) = 1/A_{d-1}^2$ on $u, v \in S^{d-1}$. As $\rho$ tends to 1, it can be shown that $P(\|U - QV\| < \varepsilon) \to 1$ for any $\varepsilon > 0$. We note that this property holds not only for density (2.1), but also for any $O(d)$-symmetric density with shape parameter $\kappa$, namely, $f(u, v) = g(u'Qv; \kappa)$. For any such model, it holds that $P(\|U - QV\| < \varepsilon)$ is monotonically increasing with respect to $\kappa$.

As is clear from the form of (2.1), $c(u, v)$ is a function of $u'Qv$, the inner product of $u$ and $Qv$. From this fact, we easily find that the density (2.1) takes maximum (minimum) values for a given $v$ at $u = Qv$ $(u = -Qv)$. The parameter $Q$ thus controls the mode of the density. It is known that an orthogonal transformation $Q$ involves two types of transformation, namely, rotation and/or reflection. In particular, when $d = 2$, these transformations can be expressed as

$$v \longmapsto \begin{pmatrix} \cos\theta & -\sin\theta \\ \sin\theta & \cos\theta \end{pmatrix} v \quad \text{and} \quad v \longmapsto \begin{pmatrix} 1 & 0 \\ 0 & -1 \end{pmatrix} v,$$

where $0 \le \theta < 2\pi$. If $\det Q = 1$, this transformation consists of only rotation. Otherwise, if $\det Q = -1$, the transformation is made up of a reflection together with a rotation.

# 3. Properties of and inference for the proposed model

## 3.1. Marginals and conditionals

One important feature of the proposed model is that it has well-known marginals and conditionals. Suppose $(U, V) \sim BS_d(\rho Q)$. The density for this random vector, (2.1), is $O(d)$- symmetric in the sense of Rivest ([1984](#)), Example 1. It then follows that the marginals of $U$ and $V$ are uniform distributions on $S^{d-1}$ with density

$$f(x) = \frac{1}{A_{d-1}}, \qquad x \in S^{d-1}.$$

Hence, model ([2.1](#)) can be viewed as a copula on $S^{d-1} \times S^{d-1}$. One difference between this special copula and the usual ones is the periodicity of its variables. As discussed in Saw ([1983](#)), Section 4, it is possible to obtain the model with rotationally symmetric marginals from the one with uniform marginals. One can also generalize the bivariate circular model, that is, $d = 2$, so that both marginals have specified distributions by using copula theory. We discuss generalizations using the above techniques and some other methods in Section [5](#).

Both conditional distributions of $U$ given $V = v$ and $V$ given $U = u$ are the exit distributions for the sphere. The terminology *exit distribution* is taken from Durrett ([1984](#)), Section 1.10, and the exit distribution on $S^{d-1}$, $\mathrm{Exit}_d(\eta)$, is of the form

$$f(x) = \frac{1}{A_{d-1}} \frac{1 - \|\eta\|^2}{\|x - \eta\|^d}, \qquad x \in S^{d-1}, \tag{3.1}$$

where $\eta \in \{\zeta \in \mathbb{R}^d; \|\zeta\| < 1\}$. This distribution is unimodal and rotationally symmetric about $x = \eta/\|\eta\|$ with the concentration being controlled by $\|\eta\|$. In particular, when $\|\eta\| = 0$, the distribution reduces to the spherical uniform. As $\|\eta\| \to 1$, it tends to a degenerate distribution concentrated at $x = \eta$. It follows that the conditionals of the model with density ([2.1](#)) are $U|(V = v) \sim \mathrm{Exit}_d(\rho Q v)$ and $V|(U = u) \sim \mathrm{Exit}_d(\rho Q' u)$.

It is worth remarking that the conditional of $W \equiv v'Q'U$ given $V = v$ has a family discussed by Leipnik ([1947](#)) and McCullagh ([1989](#)). The derivation of the density is clear from the latter paper or from Watson's result (Watson ([1983](#)), Equation 3.4.1). As in the latter paper, write $X \sim H'(\theta, \nu)$ if the density of the random variable $X$ is

$$f(x) = \frac{1 - \theta^2}{B(\nu + 1/2, 1/2)} \frac{(1 - x^2)^{\nu - 1/2}}{(1 - 2\theta x + \theta^2)^{\nu + 1}}, \qquad -1 < x < 1,$$

where $-1 < \theta < 1$ and $\nu > -1/2$. It then follows that $W|(V = v) \sim H'\{\rho, (d-2)/2\}$.

## 3.2. Some properties

Here, we investigate some of the properties of the model with density ([2.1](#)). The first is that the distribution is closed under orthogonal transformations:

$$(U, V) \sim BS_d(\rho Q) \quad \Longrightarrow \quad (Q_1 U, Q_2 V) \sim BS_d(\rho Q_1 Q Q_2'), \qquad Q_1, Q_2 \in O(d).$$

The next result can be obtained by applying a result which appears in, for example, Durrett ([1984](#)), Section 1.10.

**Theorem 1.** *Suppose that $(U, V)$ is distributed as $BS_d(\rho Q)$. Let $f$ be $C^2$ in $D$ and continuous on $\overline{D}$, where $D = \{\zeta \in \mathbb{R}^d; \|\zeta\| < 1\}$. If $f$ is harmonic, namely,*

$$\frac{\partial^2}{\partial x_1^2}f + \frac{\partial^2}{\partial x_2^2}f + \cdots + \frac{\partial^2}{\partial x_d^2}f = 0,$$

*then $E\{f(V)|U = u\} = f(\rho Q'u)$ and $E\{f(U)|V = v\} = f(\rho Qv)$.*

Using this fact, it is easy to show that $E\{f(U)\} = E\{f(V)\} = f(0)$.

The moments and correlation coefficient of the model are given by the following theorem.

**Theorem 2.** *Suppose $(U, V)$ has density (2.1). Then*

$$E(U) = E(V) = 0, \quad E(UU') = E(VV') = d^{-1}I,$$
$$E(UV') = d^{-1}\rho Q. \tag{3.2}$$

*The Johnson and Wehrly (1977) coefficient of correlation, $\rho_{JW}$, is thus*

$$\rho_{JW} \equiv \lambda^{1/2} = \rho,$$

*where $\lambda$ is the largest eigenvalue of $\Sigma_{UU}^{-1}\Sigma_{UV}\Sigma_{VV}^{-1}\Sigma_{UV}', \Sigma_{UU} = E(UU') - E(U)E(U'),$ $\Sigma_{UV} = E(UV') - E(U)E(V')$ and $\Sigma_{VV} = E(VV') - E(V)E(V')$.*

See Appendix B for the proof.

Note the simplicity of these moments and of the correlation coefficient. To our knowledge, no distributions for bivariate angular data have such a simple correlation coefficient except for the uniform distribution.

The following is useful for constructing a pivotal statistic for $(\rho, Q)$, which is discussed in Section 3.5. The result is stated in its general form as follows. The proof is also given in Appendix B.

**Theorem 3.** *Assume $(U, V)$ has $O(d)$-symmetric density, that is, $f(u, v) = g(u'Qv)$, where $Q \in O(d)$. The density of $T \equiv U'QV$ is then given by*

$$h(t) = A_{d-1} \cdot A_{d-2}g(t)(1-t^2)^{(d-3)/2}, \qquad -1 \le t \le 1.$$

From this theorem, it follows that if $(U, V) \sim BS_d(\rho Q)$, then $U'QV \sim H'\{\rho, (d-1)/2\}$. Note that this result also enables us to obtain the density for $X'QY$ of Saw's (1983) distribution immediately.

## 3.3. Random vector simulation

To generate a random vector having density (2.1), it is helpful to apply the idea which appears in, for instance, Saw (1983) and Watson (1983), Equation 2.2.1. We generalize the idea as follows.

Let $W$ be a random variable from $H'\{\rho,(d-2)/2\}$ and let $d(X;\zeta) = (I-\zeta\zeta')X/\|(I-\zeta\zeta')X\|$, where $\zeta \in S^{d-1}$ and $X$ is a random vector having a uniform distribution on $S^{d-1}$. In other words, $d(X;\zeta)$ has a uniform distribution on the $(d-1)$-sphere, $S_\perp$, in $\mathbb{R}^d$ defined by $S_\perp = \{\eta \in \mathbb{R}^d; \|\eta\| = 1, \zeta'\eta = 0\}$. The conditional of $U$ given $V = v$ can then be decomposed into

$$U|(V=v) \stackrel{\mathrm{d}}{=} WQv + (1-W^2)^{1/2}d(X;Qv).$$

Given this, generation of variates from (2.1) can be carried out using the following three steps: (i) generate a random vector $V$ which has a uniform distribution on $S$, achieved by using the method proposed by Tashiro (1977); (ii) generate $W$, which has $H'\{\rho,(d-2)/2\}$, as stated in Section 4 of McCullagh (1989); (iii) finally, a random vector $d(X;Qv)$ distributed as a uniform distribution on $S_\perp$ is obtained in a similar manner as in step (i) and one obtains a variate from the conditional of $U$ given $V = v$ as described in the preceding paragraph. The joint distribution of $(U,V)$ is then $BS_d(\rho Q)$.

## 3.4. Parameter estimation

Parameter estimation for multivariate distributions is often difficult. This is also the case for our model. However, one can discuss parameter estimation under certain conditions. Here, we consider parameter estimation based on the method of moments and maximum likelihood.

First, the method of moments estimator is constructed from (3.2). Assume that $(U_j, V_j)$ $(j = 1, \ldots, n\,(\geq 2))$ is a random sample from a distribution with density (2.1) with unknown parameters $\rho$ and $Q$. Under the condition $\operatorname{rank}(\sum_{j=1}^n U_j V_j') = d$, one can construct an estimator for the parameters based on the moment $E(UV')$. This is done by equating the theoretical and sample moments. We thus obtain

$$\hat{\rho}\hat{Q} = \frac{d}{n}\sum_{j=1}^n U_j V_j'. \tag{3.3}$$

The estimators of $\rho$ and $Q$ induced from the condition $|\det Q| = 1$ are then given by

$$\hat{\rho} = d\left|\det\left(\frac{1}{n}\sum_{j=1}^n U_j V_j'\right)\right|^{1/d} \quad \text{and} \quad \hat{Q} = \frac{d}{n\hat{\rho}}\sum_{j=1}^n U_j V_j'.$$

The estimator $\hat{\rho}\hat{Q}$ has the following properties. For the proof, see Appendix B.

**Theorem 4.** *The following hold for the estimator $\hat{\rho}\hat{Q}$ defined in (3.3):*

(i) *$\hat{\rho}\hat{Q}$ is an unbiased and consistent estimator of $\rho Q$;*

(ii) *if g is a function defined by* $g(A) = (a'_1, \ldots, a'_d)'$, *where* $A = (a_1, \ldots, a_d)$ *is a* $d \times d$ *matrix, then*

$$\sqrt{n}\{g(\rho\hat{Q}) - g(\rho Q)\} \xrightarrow{d} N(0, \Sigma) \qquad as \ n \to \infty,$$

*where* $\Sigma = (\sigma_{mn})$ *is*

$$\sigma_{mn} = \begin{cases} 1 + \dfrac{\rho^2}{d+2}\{(d-2)q_{ij}^2 - 2\}, & m = n, \\[2mm] \dfrac{\rho^2}{d+2}(dq_{kj}q_{il} - 2q_{ij}q_{kl}), & otherwise, \end{cases} \tag{3.4}$$

$q_{ij}$ *is the* $(i,j)$*th entry of* $Q$, $m = d(j-1)+i$ *and* $n = d(l-1)+k$, $1 \le i, j, k, l \le d$.

We note that although $\hat{Q}$ is an unbiased estimator of $Q$ with $|\det\hat{Q}| = 1$, it is not necessarily an orthogonal matrix.

Next, we consider maximum likelihood estimation. Let $(U_j, V_j)$ $(j = 1, \ldots, n)$ be an i.i.d. sample from $BS_d(\rho Q)$, where $Q$ is known and $\rho$ is unknown. The log-likelihood for $\rho$ is given by

$$l(\rho) = C + n\log(1 - \rho^2) - \frac{d}{2}\sum_{j=1}^{n}\log(1 - 2\rho u'_j Q v_j + \rho^2), \tag{3.5}$$

where $C$ is a constant which does not depend on $\rho$. The derivative with respect to $\rho$ is

$$\frac{\partial l}{\partial \rho} = -\frac{2n\rho}{1 - \rho^2} + d\sum_{j=1}^{n}\frac{x_j - \rho}{1 - 2\rho x_j + \rho^2},$$

where $x_j = u'_j Q v_j \in [-1, 1]$. From this expression, we find that the maximization of (3.5) with respect to $\rho$ is essentially the same as that of $H'\{\rho, (d-2)/2\}$ with respect to $\rho$.

## 3.5. Pivotal statistic

Suppose $(U, V)$ is a $BS_d(\rho Q)$ random vector. Define a random variable

$$T(\rho, Q) = \frac{1 - (U'QV)^2}{1 - 2\rho U'QV + \rho^2}.$$

It is easy to see that $0 < T(\rho, Q) < 1$ a.s. for any $\rho$ and $Q$. As shown in Theorem 3, $U'QV \sim H'\{\rho, (d-2)/2\}$. By using the results in McCullagh (1989), Section 4, and Abramowitz and Stegun (1972), equations (15.1.13) and (15.3.1), one then obtains

$$E\{T(\rho, Q)^r\} = \frac{B(r + (d-1)/2, 1/2)}{B((d-1)/2, 1/2)},$$

where $B(\cdot, \cdot)$ is a beta function. Since these moments are equal to those of a beta distribution Beta$((d-1)/2, 1/2)$, it follows that $T(\rho, Q)$ is a pivotal statistic for $(\rho, Q)$ having a Beta$((d-1)/2, 1/2)$ distribution almost surely.

# 4. Bivariate circular case

## 4.1. Transformation of random vectors and parameters

Thus far, we have considered properties of model (2.1) for the general dimensional case. In this section, we specifically discuss the bivariate circular case of the proposed model which possesses some unique properties.

Suppose $(U, V) \sim BS_2(\rho Q)$. Its density is then expressed as

$$c(u, v) = \frac{1}{4\pi^2} \frac{1 - \rho^2}{1 - 2\rho u' Q v + \rho^2}, \qquad u, v \in S^1.$$

For ease of discussion, it will be helpful to transform the random variables and parameters by taking

$$(Z_U, Z_V) = (U_1 + iU_2, V_1 + iV_2) \quad \text{and} \quad \psi = \rho e^{i\theta},$$

where $U = (U_1, U_2)'$, $V = (V_1, V_2)'$ and $\theta$ is a constant satisfying

$$Q = \begin{pmatrix} \cos\theta & -\det Q \sin\theta \\ \sin\theta & \det Q \cos\theta \end{pmatrix}, \qquad 0 \le \theta < 2\pi. \tag{4.1}$$

It then follows that $|\psi| < 1$ and $Z_U, Z_V \in \Omega$, where $\Omega = \{z \in \mathbb{C}; |z| = 1\}$. The density for $(Z_U, Z_V)$ is given by

$$c(z_u, z_v) = \frac{1}{4\pi^2} \frac{1 - |\psi|^2}{|1 - \psi z_v z_u^{-\det Q}|^2}, \qquad z_u, z_v \in \Omega. \tag{4.2}$$

If $(Z_U, Z_V)$ has density (4.2) with $\det Q = 1$, we write $(Z_U, Z_V) \sim BC_+(\psi)$. Similarly, we write $(Z_U, Z_V) \sim BC_-(\psi)$ if $(Z_U, Z_V)$ has density (4.2) with $\det Q = -1$.

Note that this transformation does not actually change the distribution. All we have done is to express the random variables and the parameters in the form of complex numbers for the sake of further investigation of the distributions.

As already stated in Section 2.2, the marginals of $Z_U$ and $Z_V$ are circular uniform, whereas both conditionals of $Z_U$ given $Z_V = z_v$ and $Z_V$ given $Z_U = z_u$ are exit distributions for the circle, that is, wrapped Cauchy distributions. For brevity, we introduce the notation $C^*(\phi)$ derived from McCullagh (1996) which denotes the wrapped Cauchy or circular Cauchy distribution with density

$$f(z) = \frac{1}{2\pi} \frac{1 - |\phi|^2}{|z - \phi|^2}, \qquad z \in \Omega; |\phi| < 1.$$

The relationship $|\phi| = \|\xi\|$ and $\arg(\phi) = \arg(\xi_1 + i\xi_2)$, where $\xi = (\xi_1, \xi_2)'$, holds between the parameters of model (3.1) and those of the density above via a transformation $Z = X_1 + iX_2$. See McCullagh (1996), Mardia and Jupp (2000), pp. 51–52, and Jammalamadaka and SenGupta (2001), pp. 45–46, for further properties of the wrapped Cauchy and circular Cauchy distribution. For model (4.2), it is easy to show that $Z_U | (Z_V = z_v) \sim C^*(\psi z_v)$ and $Z_V | (Z_U = z_u) \sim C^*(\overline{\psi} z_u)$.

## 4.2. Some properties

To investigate other properties of the model, it is useful to calculate its moments. Assume that $(Z_U, Z_V)$ has $BC_+(\psi)$. The moments for $(Z_U, Z_V)$ are then obtained, by applying Rudin (1987), Theorem 11.13, as

$$E(Z_U{}^j Z_V{}^k) = \begin{cases} \psi^j, & j = -k \geq 0, \\ \overline{\psi}^{-j}, & j = -k < 0, \\ 0, & \text{otherwise,} \end{cases} \qquad \text{for } j, k \in \mathbb{Z}. \tag{4.3}$$

Similarly, we can obtain the moments for $BC_-$. According to Fourier series expansion theory, one can recover the density from these moments if the density $f$ satisfies $f \in L^2(\Omega \times \Omega)$. See Dym and McKean (1972), Section 1.10, for details.

Using these results, the following properties are established. First, the $BC_+$ model has the multiplicative property

$$(Z_{U1}, Z_{V1}) \sim BC_+(\psi_1) \perp (Z_{U2}, Z_{V2}) \sim BC_+(\psi_2)$$

$$\implies (Z_{U1} Z_{U2}, Z_{V1} Z_{V2}) \sim BC_+(\psi_1 \psi_2). \tag{4.4}$$

Likewise, it can be shown that the $BC_-$ model also has this multiplicative property. However, the convolution of $BC_+$ and $BC_-$ is the uniform distribution, that is,

$$(Z_{U1}, Z_{V1}) \sim BC_+(\psi_1) \perp (Z_{U2}, Z_{V2}) \sim BC_-(\psi_2)$$

$$\implies (Z_{U1} Z_{U2}, Z_{V1} Z_{V2}) \sim BC_+(0).$$

In addition,

$$(Z_U, Z_V) \sim BC_\pm(\psi) \implies (Z_U{}^n, Z_V{}^n) \sim BC_\pm(\psi^n) \qquad \text{for any } n \in \mathbb{N}.$$

As $n$ tends to infinity, the distribution of $(Z_U{}^n, Z_V{}^n)$ tends to a uniform distribution on the torus.

Furthermore, model (4.2) is infinitely divisible with respect to multiplication. This can be proven as follows. Let $(Z_U, Z_V) \sim BC_\pm(\psi)$. For any positive integer $n$, the assumption that $(Z_{Uj}, Z_{Vj})$ $(j = 1, \ldots, n)$ is an i.i.d. sample from $BC_\pm(\sqrt[n]{\psi})$ then yields

$$\left( \prod_{j=1}^n Z_{Uj}, \prod_{j=1}^n Z_{Vj} \right) \stackrel{d}{=} (Z_U, Z_V). \tag{4.5}$$

### 4.3. Random vector simulation

In order to simulate a $BC_+(\psi)$ random vector, one could generate $\mathbb{R}^2$-valued Brownian motion and record the points at which the Brownian particle hits circles with radii $\rho$ and 1. However, this algorithm is somewhat inefficient because we need to simulate Brownian motion at least up to the time at which the particle hits the unit circle. Another possibility is discussed in Section 3.3, but it, too, is less efficient than the method proposed below. The focus of this subsection is therefore to discuss an algorithm for simulating $BC_+(\psi)$ variates which we conclude to be more appealing than the aforementioned methods.

To obtain the random vector, we use the fact that the marginal of $Z_U$ is circular uniform and the conditional of $Z_V$ given $Z_U = z_u$ is wrapped Cauchy, specifically, $C^*(\overline{\psi}z_u)$. For the generation of a variate from a wrapped Cauchy distribution, we apply a result from McCullagh (1996) concerning the Möbius transformation of a circular uniform random variable, namely that

$$Z \sim C^*(0) \quad \Longrightarrow \quad \frac{Z + \beta}{1 + \overline{\beta}Z} \sim C^*(\beta), \qquad |\beta| < 1. \tag{4.6}$$

An algorithm for generating $BC_+(\psi)$ random vectors then involves the following steps.

Step 1: Generate uniform $(0, 1)$ random numbers $U_1$ and $U_2$.
Step 2: Set $Z_U = \exp(2\pi \mathrm{i}U_1)$ and $Z_T = \exp(2\pi \mathrm{i}U_2)$.
Step 3: Let $Z_V = \frac{\overline{\psi}Z_U + Z_T}{1 + \psi \overline{Z_U} Z_T}$.

The joint distribution of $(Z_U, Z_V)$ is then $BC_+(\psi)$. In Step 2, $Z_U$ and $Z_T$ are independent circular uniform random variables. In Step 3, because of property (4.6), the conditional distribution of $Z_V$ given $Z_U = z_u$ is $C^*(\overline{\psi}z_u)$. It therefore follows that $(Z_U, Z_V) \sim BC_+(\psi)$.

$BC_-(\psi)$ random vectors can be simulated using a very similar approach.

### 4.4. Parameter estimation

Here, we consider parameter estimation for the $BC_+(\psi)$ model based on the method of moments and maximum likelihood. Although we discuss parameter estimation for only the $BC_+(\psi)$ here, it is possible to derive the estimates of the parameters for the $BC_-(\psi)$ model by a straightforward modification of the result below.

First, we consider method of moments estimation based on (4.3). Assume $(Z_U, Z_V)$ is a $BC_+(\psi)$ random variable. As discussed in Section 4.2, its theoretical moments are given by (4.3). Suppose $(Z_{Uj}, Z_{Vj})$ $(j = 1, \ldots, n)$ is a random sample from the $BC_+(\psi)$ distribution. The method of moments estimator is obtained by equating the theoretical and sample moments. We thus obtain

$$\hat{\psi} = \frac{1}{n} \sum_{j=1}^{n} Z_{Uj} \overline{Z_{Vj}}. \tag{4.7}$$

It follows from the weak law of large numbers that $\hat{\psi}$ is a consistent estimator of $\psi$. In addition, the central limit theorem enables us to prove asymptotic normality of the estimator, namely,

$$\sqrt{n}\left\{\begin{pmatrix} \mathrm{Re}(\hat{\psi}) \\ \mathrm{Im}(\hat{\psi}) \end{pmatrix} - \begin{pmatrix} \mathrm{Re}(\psi) \\ \mathrm{Im}(\psi) \end{pmatrix}\right\} \xrightarrow{d} N\left\{0, \frac{1}{2}(1-|\psi|^2)I\right\} \qquad \text{as } n \to \infty.$$

Although this estimator is different from the method of moments estimator (3.3), these estimators are somewhat related. Recall that the relationship (4.1) holds between $\arg(\psi)$ and $Q$. If $\det Q = 1$, $\rho Q$ can then be expressed as

$$\rho Q = \begin{pmatrix} \mathrm{Re}(\psi) & -\mathrm{Im}(\psi) \\ \mathrm{Im}(\psi) & \mathrm{Re}(\psi) \end{pmatrix}.$$

Given this relationship, an estimator of $\psi$ induced naturally from the method of moments estimator (3.3) is

$$\hat{\psi} = \mathrm{Re}\left\{\frac{1}{n}\sum_{j=1}^{n}(U_{j1}V_{j1} + U_{j2}V_{j2})\right\} + \mathrm{i}\,\mathrm{Im}\left\{\frac{1}{n}\sum_{j=1}^{n}(U_{j2}V_{j1} - U_{j1}V_{j2})\right\},$$

where $U_j = (U_{j1}, U_{j2})'$ and $V_j = (V_{j1}, V_{j2})'$. This estimator is equal to the method of moments estimator (4.7).

Second, turning to the maximum likelihood estimation, it is obvious that the maximum likelihood estimator coincides with the method of moments estimator, that is, $\hat{\psi} = Z_{U1}\overline{Z_{V1}}$, for a single observation (i.e., when $n = 1$). When $n$ is large, the estimates must be obtained numerically. Note that the likelihood function can be written as

$$L(\psi) \propto \prod_{j=1}^{n} \frac{1-|\psi|^2}{|z_{uj}\overline{z_{vj}} - \psi|^2}.$$

This expression suggests that maximum likelihood estimation for the $BC_+(\psi)$ model essentially coincides with that for the wrapped Cauchy distribution $C^*(\psi)$. We can therefore obtain estimates by applying the algorithm of Kent and Tyler (1988). Since distribution (4.2) is identifiable and the parameter space, the unit disc, is finite, consistency of the maximum likelihood estimator follows from the general theory (see, e.g., Bahadur (1971)). The Fisher information matrix for $(\mathrm{Re}(\psi), \mathrm{Im}(\psi))$, denoted by $\mathcal{I}(\mathrm{Re}(\psi), \mathrm{Im}(\psi))$, is simply expressed as

$$\mathcal{I}\{\mathrm{Re}(\psi), \mathrm{Im}(\psi)\} = \frac{2}{(1-|\psi|^2)^2}I.$$

The above can be obtained by transforming random variables into polar coordinates and using (3.616.7) of Gradshteyn and Ryzhik (1994). Therefore, the following hold for the maximum likelihood estimator:

$$\sqrt{n}\left\{\begin{pmatrix} \mathrm{Re}(\hat{\psi}) \\ \mathrm{Im}(\hat{\psi}) \end{pmatrix} - \begin{pmatrix} \mathrm{Re}(\psi) \\ \mathrm{Im}(\psi) \end{pmatrix}\right\} \xrightarrow{d} N\left\{0, \frac{1}{2}(1-|\psi|^2)^2 I\right\} \qquad \text{as } n \to \infty.$$

**Table 1.** Estimates of relative mean squared errors of the method of moments estimator (4.7) with respect to the maximum likelihood estimator

|  | $\psi = 0.1$ | $\psi = 0.3$ | $\psi = 0.5$ | $\psi = 0.7$ | $\psi = 0.9$ |
|---|---|---|---|---|---|
| $n = 10$ | 0.919 | 0.998 | 1.155 | 1.620 | 4.135 |
| $n = 20$ | 0.963 | 1.032 | 1.221 | 1.749 | 4.767 |
| $n = 30$ | 0.980 | 1.071 | 1.229 | 1.795 | 4.942 |
| $n = 50$ | 0.977 | 1.059 | 1.306 | 1.827 | 5.039 |
| $n = 100$ | 0.992 | 1.105 | 1.311 | 1.891 | 5.088 |
| $n = \infty$ | 1.010 | 1.099 | 1.333 | 1.961 | 5.263 |

## 4.5. Simulation study

In this subsection, a simulation study is carried out to compare the finite sample performance of the estimators. Here, we discuss two estimators for the parameter $\psi$ based on the method of moments (4.7) and maximum likelihood. As for the estimator (3.3), we do not discuss it here because it is expressed in matrix form and we cannot directly compare it with the other two estimators.

In our simulation study, random samples of sizes $n = 10, 20, 30, 50$ and 100 from $BC_+(\psi)$ with $\psi = 0.1, 0.3, 0.5, 0.7$ and 0.9 are generated. For each combination of $n$ and $\psi$, 2000 random samples are gathered. Random vectors from $BC_+(\psi)$ are generated by using a method introduced in Section 4.3. We employ the Mersenne Twister, which is implemented by the command `runif` in R 2.7.1, to obtain uniform random variates.

We discuss the performance of the estimators in terms of the estimates of the mean squared error. In this case, the estimate of the mean squared error is given by $\sum_{j=1}^{2000} |\hat{\psi}_j - \psi|^2/2000$, where the $\hat{\psi}_j$'s $(j = 1, \ldots, 2000)$ are the estimates for $\psi$.

Estimates of the relative mean squared errors of the method of moments estimator (4.7) with respect to the maximum likelihood estimator for some selected values of $n$ and $\psi$ are given in Table 1. In the table, the relative mean squared errors for $n = \infty$ are derived from the asymptotic variance of two estimators as $(1 - |\psi|^2)/(1 - |\psi|^2)^2 = 1/(1 - |\psi|^2)$. The comparison of the relative mean squared errors shows that the maximum likelihood estimator provides a better result in many cases, especially for large $|\psi|$ or $n$. However, the difference diminishes as $n$ or $|\psi|$ decreases. In particular, when $\psi = 0.1$ or $(\psi, n) = (0.3, 10)$, the method of moments estimator displays better performance.

An advantage of the method of moments estimator (4.7) is its simplicity in calculation. Because the estimator is expressed in closed form, it is not necessary to use any numerical algorithm to estimate the parameter. In addition, the estimator displays satisfactory performance when $|\psi|$ is small. However, the maximum likelihood estimator is considered a better estimator in most combinations of $\psi$ and $n$, as seen in Table 1. A small drawback of the maximum likelihood estimator is the complexity involved in calculating the estimator. Since the maximum likelihood estimator is not expressed in a closed form, we need to resort to a numerical method to obtain the estimate for the parameter.

# 5. Related models

## 5.1. Generalizations of the model with density (2.1)

As described in Section 2.1, the model with density (2.1) is generated using Brownian motion starting at $B_0 = 0$. In this subsection, we briefly discuss a distribution which is generated using Brownian motion starting at $B_0 = \xi(\|\xi\| < \rho)$ instead of $B_0 = 0$. We define a random vector $(U, V) = (QB_{\tau_1}/\|B_{\tau_1}\|, B_{\tau_2})$ in the same way as was used in Section 2.1, except that we incorporate the new starting point. The resulting density for $(U, V)$ is given by

$$f(u, v) = \frac{1}{A_{d-1}^2} \frac{1 - \rho^2}{(1 - 2\rho u' Q v + \rho^2)^{d/2}} \frac{\rho^2 - \|\xi\|^2}{(\rho^2 - 2\rho U' Q \xi + \|\xi\|^2)^{d/2}}, \qquad u, v \in S^{d-1}. \quad (5.1)$$

The marginals and conditional distribution of $V$ given $U = u$ are the exit distributions:

$$U \sim \text{Exit}_d(\rho^{-1} Q \xi), \qquad V \sim \text{Exit}_d(\xi) \quad \text{and} \quad V|(U = u) \sim \text{Exit}_d(\rho Q' u).$$

The conditional distribution of $U$ given $V = v$ is not of the usual form. This conditional distribution can be unimodal or bimodal and is generally skewed, except for certain special cases such as $v = \pm \xi/\|\xi\|$. It can be shown that $U$ and $V$ are independent if and only if $\rho = 0$. We note that the bivariate circular case of model (5.1) is a submodel of the distribution briefly discussed by Kato *et al.* (2008) as a model related to a circular-circular regression model.

Another generalization arises from the use of the method discussed in Saw (1983), Section 4. This method enables us to derive a distribution with prescribed rotationally symmetric marginals.

In the bivariate circular case, it might be promising to apply the Möbius transformation to each variable. Let $(\tilde{Z}_U, \tilde{Z}_V) \sim BC_+(\psi)$ and define a random vector

$$(Z_U, Z_V) = \left( \frac{\tilde{Z}_U + \alpha_1}{1 + \overline{\alpha_1} \tilde{Z}_U}, \frac{\tilde{Z}_V + \alpha_2}{1 + \overline{\alpha_2} \tilde{Z}_V} \right), \qquad |\alpha_1|, |\alpha_2| < 1. \quad (5.2)$$

Then, because of property (4.6), the marginals of $Z_U$ and $Z_V$ have wrapped Cauchy distributions $C^*(\alpha_1)$ and $C^*(\alpha_2)$, respectively. Another benefit of this extension is that its density has a simple and exact form, including the normalizing constant which does not involve any special functions.

It is also possible to transform the bivariate circular model into a distribution with specified marginals by applying Sklar's theorem in the theory of copulas. (See Nelsen (1998), Theorem 2.3.3.) For example, a bivariate distribution with von Mises or, equivalently, circular normal marginals is constructed as follows. Let $(Z_U, Z_V)$ have $BC_+(\psi)$. Assume that $(\tilde{\Theta}_U, \tilde{\Theta}_V) = (\text{Arg}(Z_U), \text{Arg}(Z_V))$, where $\text{Arg}(z)$ is the argument of $z$ taking values on $[0, 2\pi)$. Suppose that $F_j$ $(j = 1, 2)$ are distribution functions of the von Mises distributions vM$(\mu_j, \kappa_j)$, namely,

$$F_j(\theta) = \int_0^\theta \frac{1}{2\pi I_0(\kappa_j)} \exp\{\kappa_j \cos(t - \mu_j)\} \, \mathrm{d}t, \tag{5.3}$$

where $0 \le \mu_j < 2\pi, \kappa_j \ge 0$ and $I_0(\cdot)$ denotes the modified Bessel function of the first kind and of order 0. Define a distribution by a random vector

$$(\Theta_U, \Theta_V) = \left(F_1^{-1}\left(\frac{\tilde{\Theta}_U}{2\pi}\right), F_2^{-1}\left(\frac{\tilde{\Theta}_V}{2\pi}\right)\right).$$

The density for this random vector is of the form

$$f(\theta_u, \theta_v) = \frac{1 - |\psi|^2}{4\pi^2 I_0(\kappa_1) I_0(\kappa_2)} \exp\{\kappa_1 \cos(\theta_u - \mu_1) + \kappa_2 \cos(\theta_v - \mu_2)\}$$

$$\times [1 + |\psi|^2 - 2|\psi| \cos[2\pi\{F_1(\theta_u) - F_2(\theta_v)\} - \arg(\psi)]]^{-1}, \tag{5.4}$$

$$0 \le \theta_u, \theta_v < 2\pi.$$

It follows from Sklar's theorem that the marginals of $\Theta_U$ and $\Theta_V$ are the von Mises, vM$(\mu_1, \kappa_1)$ and vM$(\mu_2, \kappa_2)$, respectively. As is clear from the derivation, the distribution reduces to $BC_+(\psi)$ when $\kappa_1 = \kappa_2 = 0$.

Figure 1 plots some contour plots of density (5.4) for fixed values of $\mu_2, \kappa_2$ and $\arg(\psi)$ and some selected values of $\mu_1, \kappa_1$ and $|\psi|$. The comparison between Figure 1(a)–(c) suggests that $|\psi|$ influences dependence between $\Theta_U$ and $\Theta_V$, and this was mathematically validated in Theorem 2. Figures 1(a) and 1(d) imply that the concentration of the marginal distributions is controlled by $\kappa_1$ or $\kappa_2$. As seen in Figures 1(e) and 1(f), the distribution can be skewed when $\mu_1 \ne \mu_2$.

It is important to decide on conditions for the independence of two variables for a bivariate distribution. The following result provides the parameter configuration which yields independence for the bivariate family including models (5.2) and (5.4).

**Theorem 5.** *Let $(Z_U, Z_V) \sim BC_+(\psi)$. Suppose $g_1$ and $g_2$ are one-to-one mappings defined on $\Omega$ and differentiable on $\Omega$. Then $g_1(Z_U)$ and $g_2(Z_V)$ are independent if and only if $\psi = 0$.*

**Proof.** Let $(\tilde{Z}_U, \tilde{Z}_V) = (g_1(Z_U), g_2(Z_V))$. The density of $(\tilde{Z}_U, \tilde{Z}_V)$ is then calculated as

$$f(\tilde{z}_u, \tilde{z}_v) = \frac{1}{4\pi^2} \frac{1 - |\psi|^2}{|1 - \psi g_2^{-1}(\tilde{z}_v)\overline{g_1^{-1}(\tilde{z}_u)}|^2} \left|\frac{\partial g_1^{-1}(\tilde{z}_u)}{\partial \tilde{z}_u}\right| \left|\frac{\partial g_2^{-1}(\tilde{z}_v)}{\partial \tilde{z}_v}\right|, \qquad \tilde{z}_u, \tilde{z}_v \in \Omega.$$

Given this, the necessary and sufficient condition for independence is that $\psi g_2^{-1}(\tilde{z}_v)\overline{g_1^{-1}(\tilde{z}_u)}$ is a constant or a function which either depends only on $\tilde{z}_u$ or only on $\tilde{z}_v$. Since neither $g_1$ nor $g_2$ is a constant function, this condition holds if and only if $\psi = 0$. $\qquad\square$

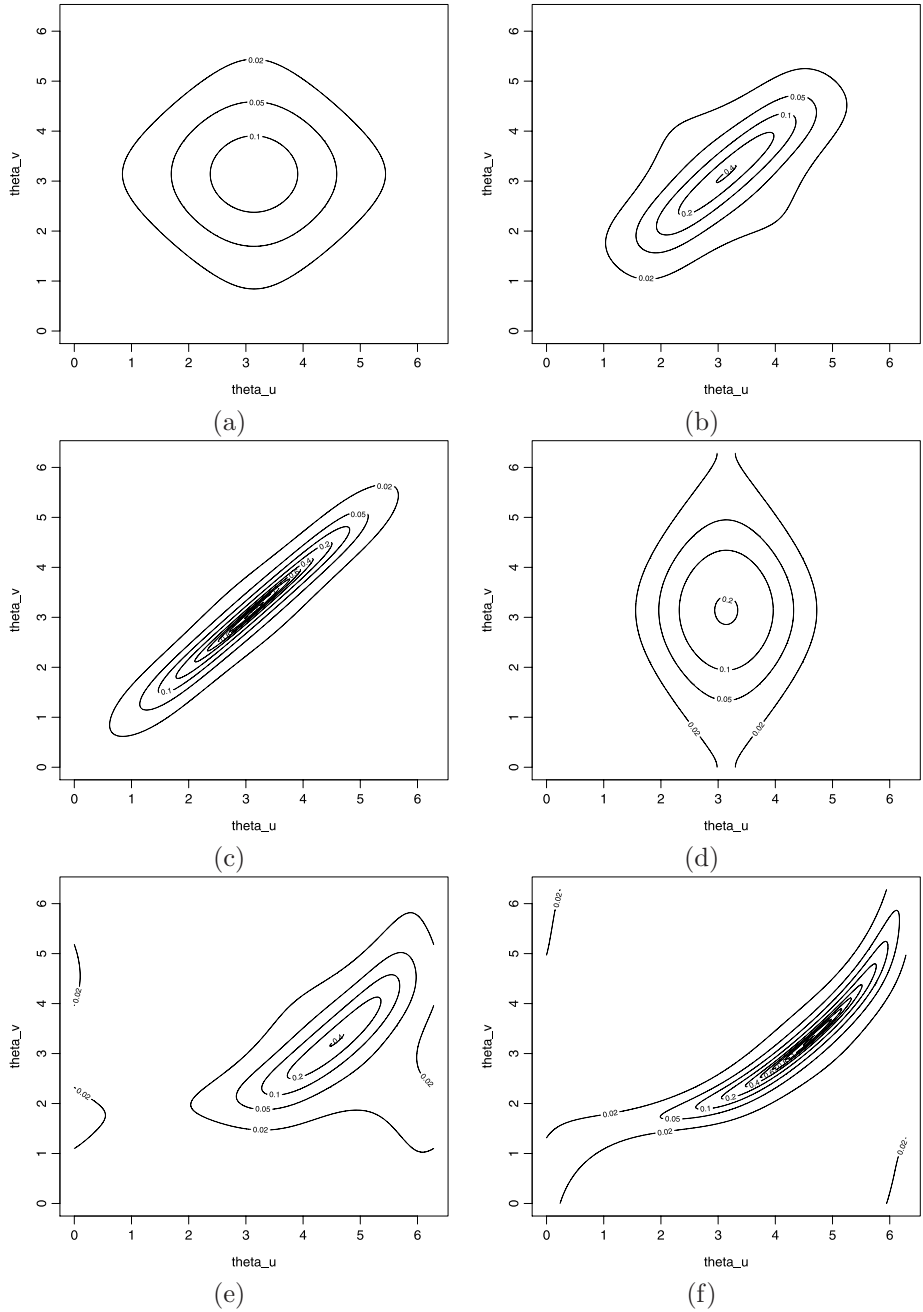

**Figure 1.** Contour plots of density (5.4) with $\mu_2 = \pi, \kappa_2 = 1.16, \arg(\psi) = 0$ and (a) $\mu_1 = \pi, \kappa_1 = 1.16, |\psi| = 0$; (b) $\mu_1 = \pi, \kappa_1 = 1.16, |\psi| = 0.5$; (c) $\mu_1 = \pi, \kappa_1 = 1.16, |\psi| = 0.8$; (d) $\mu_1 = \pi, \kappa_1 = 2.32, |\psi| = 0$; (e) $\mu_1 = 3\pi/2, \kappa_1 = 1.16; |\psi| = 0.5$; (f) $\mu_1 = 3\pi/2, \kappa_1 = 1.16, |\psi| = 0.8$.

## 5.2.  Comparison with existing bivariate circular distributions

Model (5.4) has some relation to models discussed by SenGupta (2004) and Shieh and Johnson (2005). A bivariate circular family related to the von Mises has been considered by SenGupta (2004). It has the density

$$f(\theta_u, \theta_v) \propto \exp\left\{ (1, \cos\theta_u, \sin\theta_u) \begin{pmatrix} m_{11} & m_{12} & m_{13} \\ m_{21} & m_{22} & m_{23} \\ m_{31} & m_{32} & m_{33} \end{pmatrix} \begin{pmatrix} 1 \\ \cos\theta_v \\ \sin\theta_v \end{pmatrix} \right\},$$
$$0 \leq \theta_u, \theta_v < 2\pi, \tag{5.5}$$

where $m_{jk}$ $(1 \leq j, k \leq 3)$ are the parameters, with $m_{11}$, a function of the other $m_{jk}$'s, being the normalizing constant. This model has the property that both conditionals follow the von Mises distributions, a property which our model does not have. On the other hand, our model has von Mises marginals, while model (5.5) generally does not. This difference comes from the derivations of the models. Model (5.5) is constructed by means of the conditional specification, while our model is obtained by transforming a distribution, which is generated by Brownian motion, via copula theory, so that both marginals have the von Mises distributions.

Shieh and Johnson (2005) presented a bivariate circular distribution which is called the *bivariate von Mises distribution* in their paper. The density of their model is of the form

$$f(\theta_u, \theta_v) = \frac{1}{4\pi^2 \prod_{j=1}^{3} I_0(\kappa_j)}$$
$$\times \exp[\kappa_1 \cos(\theta_u - \mu_1) + \kappa_2 \cos(\theta_v - \mu_2)$$
$$+ \kappa_3 \cos[2\pi\{F_1(\theta_u) - F_2(\theta_v)\} - \mu_3]], \qquad 0 \leq \theta_u, \theta_v < 2\pi, \tag{5.6}$$

where $0 \leq \mu_j < 2\pi, \kappa_j \geq 0, j = 1, 2, 3$, and $F_k(\cdot)$ are the distribution functions of the von Mises $\mathrm{vM}(\mu_k, \kappa_k), k = 1, 2$, as defined in (5.3). A property common to their model and ours is that both models have the von Mises marginals and belong to a general class of distributions presented by Wehrly and Johnson (1980). One difference is the conditional distributions of the model and this distinction could make a difference in fit, as we see in the next subsection.

## 5.3.  Application

***Example 1.*** As an illustrative example in which one of our models is utilized, we consider a data set of pairs of wind directions measured at a weather station in Texas. The data set is a part of a larger data set which is taken from a website <http://data.eol.ucar.edu/codiac/dss/id=85.034>. The original data set contains hourly resolution surface meteorological data from the Texas Natural Resources Conservation Commission Air Quality Monitoring Network. (These data are provided by NCAR/EOL with the support of the National Science Foundation.) Among this data set, we focus on 30 pairs of wind directions at 6 a.m. ($\theta_u$) and 7 a.m. ($\theta_v$) measured each day at a weather station, which is denoted `C28_1` in the data set, from June 1 to June 30, 2003.

**Table 2.** The maximized log-likelihood, AIC and BIC values of the proposed model (5.4) and two existing models (5.5) and (5.6) fitted to wind directions at 6 a.m. and 7 a.m.

| Model | $\log L$ | AIC | BIC |
|-------|----------|------|------|
| (5.4) | $-65.9$  | 143.8 | 152.2 |
| (5.5) | $-68.2$  | 152.4 | 163.6 |
| (5.6) | $-70.9$  | 153.8 | 162.2 |

Figure 2(a) shows a planer of the dataset. This frame suggests that there is dependence between wind directions at 6 a.m. and 7 a.m. We fit models (5.4), (5.5) and (5.6) to the data set based on maximum likelihood estimation. To estimate the parameters numerically, we employ the PORT routine which is an optimization method carried out using `nlminb` in `R 2.7.1`.

Table 2 shows the maximum log-likelihood, AIC and BIC values of the fitted models. According to the AIC and BIC criteria, our model (5.4) provides the best fit of all. Model (5.5) is the second best model judging from the AIC criterion, while the BIC criterion indicates that model (5.6) is the second best. The estimated parameters of model (5.4) are given by $\hat{\mu}_1 = 1.89, \hat{\mu}_2 = 2.01, \hat{\kappa}_1 = 1.03, \hat{\kappa}_2 = 1.19, \arg(\hat{\psi}) = 6.24$ and $|\hat{\psi}| = 0.75$. The fitted density of model (5.4) is displayed in Figure 2(b) which seems to show a satisfactory fit of the model to the data set. Since the parameter $|\hat{\psi}|$, which controls the dependence between two circular variables, is fairly large, it seems that the wind directions at 6 a.m. and 7 a.m. are strongly associated. Also, note that the argument of $\hat{\psi}$ is close to zero, implying that the mean direction of the wind directions at 6 a.m. is close to that at 7 a.m. In other words, the mean direction of $F_1(\Theta_U) - F_2(\Theta_V)$ is nearly zero. These results correspond to our intuition that wind directions usually do not change dramatically within one hour.

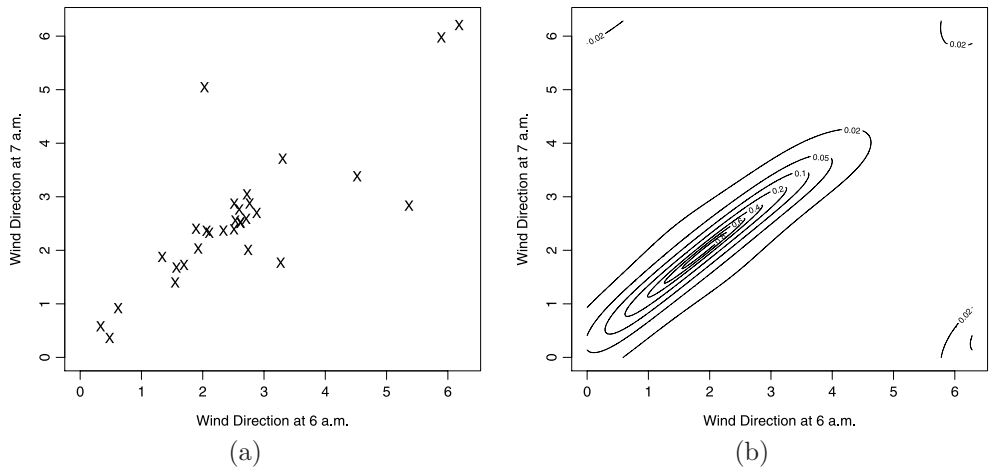

**Figure 2.** (a) Planar plot of the wind directions at 6 a.m. and 7 a.m.; (b) contour plot of the density for model (5.4) fitted to the data.

**Table 3.** The maximized log-likelihood, AIC and BIC values of the proposed model (5.4) and two existing models (5.5) and (5.6) fitted to wind directions at 6 a.m. and 6 p.m.

| Model | $\log L$ | AIC | BIC |
|-------|----------|-----|-----|
| (5.4) | −89.8 | 191.6 | 200.0 |
| (5.5) | −82.0 | 180.0 | 191.2 |
| (5.6) | −89.9 | 191.8 | 200.2 |

***Example 2.*** The second example concerns another 30 pairs of wind directions measured at the same weather station as in Example 1. This time we focus on wind directions at 6 a.m. ($\theta_u$) and 6 p.m. ($\theta_v$), observed from June 1 to June 30, 2003. Figure 3(a) shows a planar plot of the data set; it seems that the dependence structure between the two variables is not as clear as in the previous example. Again, we use models (5.4), (5.5) and (5.6) to model the data set.

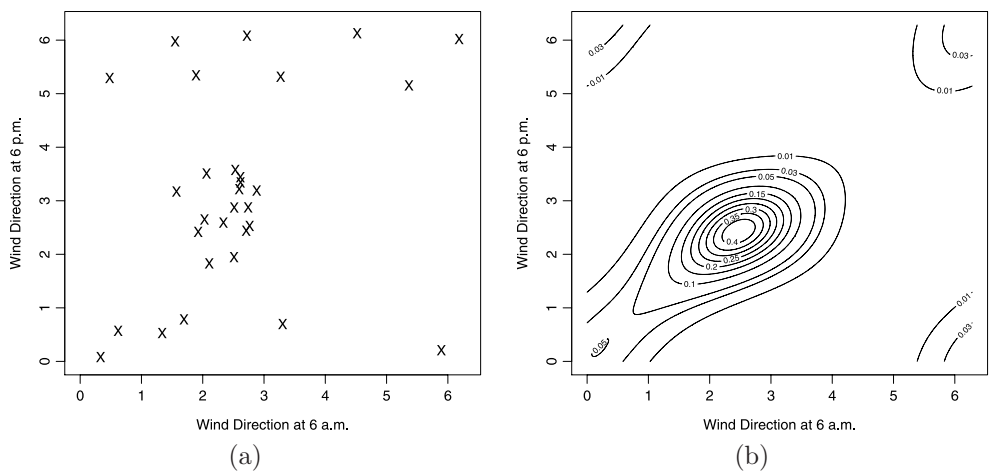

(a)                                    (b)

**Figure 3.** (a) Planar plot of the wind directions at 6 a.m. and 6 p.m.; (b) contour plot of the density for model (5.5) fitted to the data.

The maximum log-likelihood, AIC and BIC values of the three fitted models are given in Table 3. From AIC and BIC criteria, model (5.5) is the best of all. The fitted density of model (5.5) is displayed in Figure 3(b). The proposed model is the second best, but there is no significant difference in fit from Shieh and Johnson's model (5.6). The maximum likelihood estimates of the parameters of our model are given by $\hat{\mu}_1 = 2.29, \hat{\mu}_2 = 1.43, \hat{\kappa}_1 = 1.33, \hat{\kappa}_2 = 0.222, \arg(\hat{\psi}) = 0.144$ and $|\hat{\psi}| = 0.544$. Note that the estimate $|\hat{\psi}|$ in this example is smaller than that in the previous one. This suggests that there is less association between wind directions at 6 a.m. and 6 p.m. than between those at 6 a.m. and 7 a.m. However the estimate of $\arg(\psi)$ in this example is also close to zero, meaning that the mean direction of the wind directions did not make a big change although twelve hours have passed since the first observation at 6 a.m.

In Example 1, which deals with a data set with clear dependence structure, the AIC and BIC criteria suggest that our model (5.4) is the best. The example implies that the proposed model (5.4) is suitable to fit bivariate circular data if $F_1(\Theta_U) - F_2(\Theta_V)$ has a unimodal and symmetric shape. On the other hand, model (5.5) is recommended for a data set which shows a different kind of association between variables, as seen in Example 2. One advantage of model (5.5) is that it is a flexible eight-parameter model having von Mises conditionals and seems to have more potential to fit various kinds of bivariate circular data because of its flexibility. Model (5.6) has some properties common to our model; for example, the marginals of both distributions are the von Mises. However, these models have different conditionals and this difference can produce considerable distinction in fit, as demonstrated in Example 1.

# 6. Related distributions on $\mathbb{R}^2$ and on the cylinder

In previous sections, we have dealt with distributions for two directional observations. In this subsection, we provide models for two other manifolds, namely, $\mathbb{R}^2$ and the cylinder.

By applying bilinear fractional transformations to model (4.2), a distribution on $\mathbb{R}^2$ is constructed. Let $(Z_U, Z_V)$ be distributed as $BC_-(\psi)$. Define a random vector $(X, Y)$ as

$$X = \mathrm{i}\frac{1 - Z_U}{1 + Z_U} \quad \text{and} \quad Y = \mathrm{i}\frac{1 - Z_V}{1 + Z_V}.$$

Clearly, $(X, Y)$ takes values in $\mathbb{R}^2$. It is straightforward to show that the joint density for $(X, Y)$ is

$$f(x, y) = \frac{1}{\pi^2} \frac{\mathrm{Im}(\theta)}{|x + y + \theta(1 - xy)|^2}, \qquad x, y \in \mathbb{R}, \tag{6.1}$$

where $\theta = \mathrm{i}(1 - \psi)/(1 + \psi)$. Since $|\psi| < 1$, it is evident that $\mathrm{Im}(\theta) > 0$.

This model has the following properties:

$$X \sim C(i), \qquad Y \sim C(i),$$

$$X|(Y = y) \sim C\left(\frac{\theta + y}{1 - \theta y}\right), \qquad Y|(X = x) \sim C\left(\frac{\theta + x}{1 - \theta x}\right),$$

where the $C(\phi)$ notation is derived from McCullagh ([1992](#)) and denotes the Cauchy distribution on the real line with location parameter $\mathrm{Re}(\phi)$ and scale parameter $\mathrm{Im}(\phi)$. Thus, the marginals and conditionals are members of the real Cauchy family. Further properties of model ([6.1](#)) are derived using the inverse transformations $Z_U = (1+\mathrm{i}X)/(1-\mathrm{i}X)$ and $Z_V = (1+\mathrm{i}Y)/(1-\mathrm{i}Y)$ which map the real line onto the unit circle in the complex plane.

A related distribution on the cylinder $\Omega \times \mathbb{R}$ is obtained in a similar fashion. Let $(Z_U, Z_V)$ be $BC_+(\psi)$ distributions. Define a random vector

$$(Z_\Theta, X) = \left( Z_U, \mathrm{i}\frac{1 - Z_V}{1 + Z_V} \right).\tag{6.2}$$

The marginals and conditionals of $(Z_\Theta, X)$ are then

$$Z_\Theta \sim C^*(0), \qquad X \sim C(i),$$

$$Z_\Theta | (X = x) \sim C^* \left( \frac{1 + \mathrm{i}x}{1 - \mathrm{i}x} \psi \right), \qquad X | (Z_\Theta = z_\theta) \sim C \left( -\mathrm{i}\frac{1 - \overline{z_\theta}\psi}{1 + \overline{z_\theta}\psi} \right).$$

Thus, the marginals are circular uniform and standard Cauchy, while the conditionals are the wrapped Cauchy and linear Cauchy distributions, respectively.

In a manner similar to that in Section [5.1](#), one can transform the model $BS_+(\psi)$ into a family of cylindrical distributions having prescribed marginals as follows. Let $(Z_U, Z_V) \sim BS_+(\psi)$ and express these variables in terms of radians, that is, $(\Theta_U, \Theta_V) = (\mathrm{Arg}(Z_U), \mathrm{Arg}(Z_V))$. Suppose that $F_\Theta$ and $F_X$ are distribution functions of any circular and linear distributions, respectively, and are strictly increasing. A random vector defined by $(\Theta, X) = (F_\Theta^{-1}(\Theta_U/(2\pi)), F_X^{-1}(\Theta_V/(2\pi)))$ then follows a distribution on the cylinder which has marginals with distribution functions $F_\Theta$ and $F_X$. For instance, if we assume that $F_\Theta$ and $F_X$ are distribution functions of the von Mises and the normal distribution, respectively, we can construct a distribution with von Mises and normal marginals.

A straightforward modification of Theorem [5](#) yields parameter configuration for independence between two variables for the distributions presented in this subsection. For example, if a random vector $(Z_\Theta, X)$ is defined as in ([6.2](#)), then $Z_\Theta$ and $V$ are independent if and only if $\psi = 0$.

## Appendix A: Derivation of density ([2.1](#))

Let $c(u, v)$ be the joint density of $(U, V) = (QB_{\tau_1}/\|B_{\tau_1}\|, B_{\tau_2})$, which is defined in the same way as in Section [2.1](#). Note that if the density for $(U, V)$ exists, it can be expressed as

$$c(u, v) = f_U(u) g_{V|U}(v|u), \qquad u, v \in S^{d-1},$$

where $f_U$ is a density for the marginal of $U$ and $g_{V|U}$ that for the conditional of $V$ given $U = u$. Clearly, the marginal of $U$ is distributed as the uniform distribution and thus $f_U(u) = 1/A_{d-1}$. Because of the Markov property of Brownian motion, the conditional of $V$ given $U = u$ is essentially equivalent to the exit distribution for the sphere generated by Brownian motion starting at $B_0 = \rho Q'u$ (see Durrett (1984), Section 1.10). The density for the exit distribution for the sphere is known to be

$$g_{V|U}(v|u) = \frac{1}{A_{d-1}} \frac{1 - \rho^2}{\|v - \rho Q'u\|^d}, \qquad v \in S^{d-1}.$$

We thus obtain the density (2.1).

Density (5.1) is obtained by a straightforward modification of the above.

## Appendix B: Proofs of Theorems 2–4

**Proof of Theorem 2.** Since the marginals of $U$ and $V$ are uniformly distributed on the sphere, it is evident that $E(U) = E(V) = 0$ and $E(UU') = E(VV') = d^{-1}I$.

We show that $E(UV') = d^{-1}\rho I$. Because model (2.1) is $O(d)$-symmetric in the sense of Rivest (1988), calculation of $E(UV')$ is simplified by applying Proposition 1 of his paper to

$$E(UV') = \text{diag}\{E(R_j S_j)\}Q,$$

where $(R, S) \sim BS_d(\rho I), R = (R_1, \ldots, R_d)', S = (S_1, \ldots, S_d)'$. Consider the integral

$$E(R_1 S_1) = \int_{S^{d-1} \times S^{d-1}} r_1 s_1 c(r, s) \, dr \, ds = \int_{S^{d-1}} \frac{r_1}{A_{d-1}} \int_{S^{d-1}} \frac{s_1}{A_{d-1}} \frac{1 - \rho^2}{\|s - \rho r\|^d} \, ds \, dr.$$

Transforming $S$ into $\tilde{S} = PS$, where $P$ is a $d \times d$ orthogonal matrix such that $P = (r, p_2, \ldots, p_d)', p_j = (p_{j1}, \ldots, p_{jd})' \in \mathbb{R}^d$, we have

$$\int_{S^{d-1}} \frac{r_1}{A_{d-1}} \int_{S^{d-1}} \frac{s_1}{A_{d-1}} \frac{1 - \rho^2}{\|s - \rho r\|^d} \, ds \, dr$$

$$= \int_{S^{d-1}} \frac{r_1}{A_{d-1}} \int_{S^{d-1}} \frac{r_1 \tilde{s}_1 + \sum_{j=2}^d p_{j1} \tilde{s}_j}{A_{d-1}} \frac{1 - \rho^2}{(1 - 2\rho \tilde{s}_1 + \rho^2)^{d/2}} \, d\tilde{s} \, dr$$

$$= \int_{S^{d-1}} \frac{\tilde{s}_1}{dA_{d-1}} \frac{1 - \rho^2}{(1 - 2\rho \tilde{s}_1 + \rho^2)^{d/2}} \, d\tilde{s}.$$

The last equality follows from $E(R) = 0$ and $E(R_1^2) = d^{-1}$.

Then, from the fact that $X \sim H'(\theta, \nu)$ implies $E(X) = \theta$ (McCullagh (1989)), the above equation can be expressed as

$$\int_{S^{d-1}} \frac{\tilde{s}_1}{dA_{d-1}} \frac{1 - \rho^2}{(1 - 2\rho \tilde{s}_1 + \rho^2)^{d/2}} \, d\tilde{s}$$

$$= \frac{1-\rho^2}{dA_{d-1}} \frac{2\pi^{(d-1)/2}}{\Gamma\{(d-1)/2\}} \int_0^\pi \frac{\cos\theta \sin^{d-2}\theta}{(1-2\rho\cos\theta+\rho^2)^{d/2}} \, \mathrm{d}\theta$$

$$= \frac{1-\rho^2}{dB\{(d-1)/2, 1/2\}} \int_{-1}^1 \frac{t(1-t^2)^{(d-3)/2}}{(1-2\rho t+\rho^2)^{d/2}} \, \mathrm{d}t$$

$$= \frac{\rho}{d}.$$

The other elements, $E(R_j S_j)$ $(2 \le j \le d)$, are calculated in a similar way. $\qquad \square$

**Proof of Theorem 3.** The distribution function of $T$, say $H$, is given by

$$H(t) = P(T \le t) = E_V\{P(U'Qv \le t|V=v)\}$$
$$= E_{\tilde{V}}\{P(U'\tilde{v} \le t|\tilde{V} = \tilde{v})\},$$

where $\tilde{V} = QV$. Define $\tilde{U} = PU$, where $P \in O(d)$ such that $P = (\tilde{v}, p_2, \ldots, p_d)', p_j \in \mathbb{R}^d$, and one obtains

$$E_{\tilde{V}}\{P(U'\tilde{v} \le t|\tilde{V} = \tilde{v})\} = \int_{S^{d-1}} \int_{\substack{\tilde{u}_1 \le t \\ \tilde{u} \in S^{d-1}}} g(\tilde{u}_1) \, \mathrm{d}\tilde{u} \, \mathrm{d}\tilde{v}$$

$$= A_{d-1} \frac{2\pi^{(d-1)/2}}{\Gamma\{(d-1)/2\}} \int_{\substack{\cos\theta \le t \\ 0 \le \theta < \pi}} g(\cos\theta) \sin^{d-2}\theta \, \mathrm{d}\theta$$

$$= A_{d-1} \cdot A_{d-2} \int_{-1}^t g(x)(1-x^2)^{(d-3)/2} \, \mathrm{d}x.$$

Thus,

$$h(t) = A_{d-1} \cdot A_{d-2} g(t)(1-t^2)^{(d-3)/2}, \qquad -1 \le t \le 1. \qquad \square$$

**Proof of Theorem 4.** It is clear from (3.2) that $\hat{\rho}\hat{Q}$ is an unbiased estimator of $\rho Q$. Consistency of the estimator follows from the weak law of large numbers. The use of the central limit theorem enables us to prove the asymptotic normality. Here, we show that the variance–covariance matrix of $g(\hat{\rho}\hat{Q})$ is given by a matrix with entries (3.4). Suppose that $(U, V) \sim BS_d(\rho Q)$, $U = (U_1, \ldots, U_d)'$ and $V = (V_1, \ldots, V_d)'$. To calculate the variance–covariance matrix, we first consider $d^2 E(U_i V_j U_k V_l)$ for $1 \le i, j, k, l \le d$. It can be expressed as

$$d^2 E(U_i V_j U_k V_l) = d^2 E_U\{U_i U_k E_{V|U}(V_j V_l|U=u)\}$$
$$= d^2 E_U[U_i U_k E_{\tilde{V}|U}\{b_j' \operatorname{diag}(\tilde{V}_j^2) b_l|U=u\}],$$

where $\tilde{V} = (\tilde{V}_1, \ldots, \tilde{V}_d)' = PQV$, $b_j$ is the $j$th column of $PQ$ and $P$ is any $d \times d$ orthogonal matrix such that the first row of $P$ is $u'$. McCullagh (1989) showed that if $X \sim H'(\theta, \nu)$, then $E(X^2) = \{1 + (2\nu+1)\theta^2\}/\{2(\nu+1)\}$, and that $\int_{-1}^{1}(1-x^2)^{\nu-1/2}/\{(1 - 2\theta x + \theta^2)B(\nu + \frac{1}{2}, \frac{1}{2})\}\,dx = 1$. Using these results, we have

$$E_{\tilde{V}|U}\{\mathrm{diag}(\tilde{V}_j^2)|U = u\} = \frac{1 - \rho^2}{d}I + \rho^2\Delta_1,$$

where $\Delta_k = (\delta_{ij})$ is a $d \times d$ matrix whose entries are given by $\delta_{ij} = 1$ for $(i, j) = (k, k)$ and $\delta_{ij} = 0$ otherwise. Therefore,

$$d^2 E(U_i V_j U_k V_l) = d^2 E_U\left\{U_i U_k q_j'\left(\frac{1 - \rho^2}{d}I + \rho^2 UU'\right)q_l\right\},$$

where $q_j$ is the $j$th column of $Q$. If $i = k$ and $j = l$, we have

$$d^2 E(U_i V_j U_k V_l) = (1 - \rho^2) + \frac{d\rho^2}{d+2}(1 + 2q_{ij}^2).$$

The above follows from the fact that $E_U(UU') = d^{-1}I$ and $E_U(U_i UU') = (2\Delta_i + I)/\{d(d+2)\}$. If $i \neq k$ or $j \neq l$, it is fairly easy to show that

$$d^2 E(U_i V_j U_k V_l) = d^2 \rho^2 q_j' E_U(U_i U_k UU')q_l = \frac{d\rho^2}{d+2}(q_{ij}q_{kl} + q_{kj}q_{il}).$$

Thus, we obtain $d^2 E(U_i V_j U_k V_l)$ for any $1 \leq i, j, k, l \leq d$. On the other hand, it follows immediately from (3.2) that $dE(U_i V_j) = \rho q_{ij}$. Summarizing these results, we obtain Theorem 4. $\qquad\square$

# Acknowledgements

I am grateful to the Associate Editor and to a referee for the comments which improved the quality of the paper. My thanks also go to Professor Kunio Shimizu and Dr. Arthur Pewsey for their helpful comments and suggestions. Financial support for this research was provided, in part, by the Ministry of Education, Culture, Sport, Science and Technology in Japan via a grant in aid of the 21st Century Center of Excellence for Integrative Mathematical Sciences: Progress in Mathematics Motivated by Social and Natural Sciences.